\documentclass[11pt]{amsart}
\usepackage{amsfonts,amssymb}
\usepackage[mathscr]{eucal}
\usepackage{amsmath, amsthm}
\usepackage{mathrsfs}
\usepackage{amsbsy}
\usepackage{dsfont}
\usepackage{cmbright}
\usepackage{bbm}
\makeindex
\input xypic
\xyoption{all}
\begin{document}
\baselineskip = 5mm
%%%%%%%%%%%%% arrows %%%%%%%%%%%%%%%%%%%%%%%%%%%%%%%%%%%%%%%%%%%%%%
\newcommand \lra {\longrightarrow}
\newcommand \hra {\hookrightarrow}
%%%%%%%%%%%%% numbers %%%%%%%%%%%%%%%%%%%%%%%%%%%%%%%
\newcommand \ZZ {{\mathbb Z}} % integers or the free group over a finite set
\newcommand \NN {{\mathbb N}} % natural numbers
\newcommand \FF {{\mathbb F}} % finite fields
\newcommand \QQ {{\mathbb Q}} % rational numbers
\newcommand \RR {{\mathbb R}} % real numbers
\newcommand \CC {{\mathbb C}} % complex numbers
%%%%%%%%%%%%%% symmetric groups %%%%%%%%%%%%%%%%%%%%%%%%%
\newcommand \sgr {{\Sigma }} % symmetric group
%%%%%%%%%%%%%%%%% categories %%%%%%%%%%%%%%%%%%%%%%%%%%%%%
\newcommand \aA {{\bf {A}}}
\newcommand \aC {{\mathscr {C}}}
\newcommand \aT {{\bf {T}}}
%%%%%%%%%%%%% symbols for general categories %%%%%%%%%%%%%%%%%%%%%%
\newcommand \im {{\rm im}}
\newcommand \Hom {{\rm Hom}}
\newcommand \colim {{\rm {colim}\, }} % colimit
\newcommand \iHom {{\underline {\rm Hom}}} % simplicial resolution
\newcommand \End {{\rm {End}}}
\newcommand \Ob {{\rm {Ob}}}
\newcommand \coker {{\rm {coker}}}
\newcommand \id {{\rm {1}}}
\newcommand \cone {{\rm {cone}}}
\newcommand \de {{\vartriangle }}
%%%%%%%%%%%%%% well recognized categories %%%%%%%%%%%%%%%%%%%%%%%%%%%%%%%%%%%%%%%%%
\newcommand \DM {{\bf {DM}}}
\newcommand \DMG {\bf {DM}'} % triangulated category of geometrical motives
\newcommand \DMGA {{{\bf {DM}}^{\otimes }_{\leq 1}}} % abelian part of \DMG
\newcommand \CHM {{\bf {M}}}
\newcommand \SH {{\bf SH}}
\newcommand \Spt {{Spt}}
\newcommand \uno {{\mathbbm 1}} % unite in a tensor category
\newcommand \Le {{\mathbbm L}} % Lefschetz object
\newcommand \Rings {{\sf Rings}}
\newcommand \modules {{{mod}}}
%%%%%%%%%%%%% algebraic geometry symbols %%%%%%%%%%%%%%%%%%%%%%%%%%
\newcommand \AF {{\mathbb A}} % affine space
\newcommand \PR {{\mathbb P}} % projective space
\newcommand \Spec {{\rm {Spec}}}
\newcommand \SP {{\mathscr {S}\! \! \mathscr {P}}} % smooth projective schemes
\newcommand \Sch {{Sch}} % schemes
\newcommand \Var {{Var}}
\newcommand \Alb {{\rm {Alb}}}
\newcommand \Sym {{\rm {Sym}}}
\newcommand \Gm {{{\mathbb G}_m}}
%%%%%%%%%%%%% some common symbols %%%%%%%%%%%%%%%%%%%%%%%%%%%%%%%
\newcommand \cha {{\rm {char}}} % characteristic
\newcommand \tr {{\rm {tr}}} % trace
\newcommand \res {{\rm {res}}} % restriction
%%%%%%%%%%%%% theorems, lemmas, etc %%%%%%%%%%%%%%%%%%%%%%%%%%%%%%%
\newtheorem{theorem}{Theorem}
\newtheorem{lemma}[theorem]{Lemma}
\newtheorem{sublemma}[theorem]{Sublemma}
\newtheorem{corollary}[theorem]{Corollary}
\newtheorem{example}[theorem]{Example}
\newtheorem{exercise}[theorem]{Exersize}
\newtheorem{proposition}[theorem]{Proposition}
\newtheorem{remark}[theorem]{Remark}
\newtheorem{notation}[theorem]{Notation}
\newtheorem{definition}[theorem]{Definition}
\newtheorem{conjecture}[theorem]{Conjecture}
%%%%%%%%%%%%%%%%%%%%%%%%%%%%%%%%%%%%%%%%%%%%%%%%%%%%%%%%%%%%%%%%%%%%%%%%%%%
\newenvironment{pf}{\par\noindent{\em Proof}.}{\hfill\framebox(6,6)
\par\medskip}
%%%%%%%%%%%%%%%%%%%%%%%%%%%%%%%%%%%%%%%%%%%%%%%%%%%%%%%%%%%%%%%%%%%%%%%%%%%
\title{\bf Zeta functions in triangulated categories}
\author{V. Guletski\v \i
%\thanks{The authors acknowledge }
}

\date{May 2, 2006}

\begin{abstract}
\noindent We prove $2$-out-of-$3$ property for rationality of
motivic zeta function in distinguished triangles in Voevodsky's
category $\DM $. As an application, we show rationality of motivic
zeta functions for all varieties whose motives are in the thick
triangulated monoidal subcategory generated by motives of
quasi-projective curves in $\DM $. Joint with a result of
P.O'Sullivan it also gives an example of a variety whose motive is
not finite-dimensional while the motivic zeta function is rational.
\end{abstract}


%16E20 Grothendieck groups, $K$-theory, etc. [See also 18F30, 19Axx,
%19D50]
%
%18D10 Monoidal categories (= multiplicative categories), symmetric
%monoidal categories, braided categories [See also 19D23]
%
%19E15 Algebraic cycles and motivic cohomology [See also 14C25,
%14C35]
%
%55U35 Abstract and axiomatic homotopy theory


\maketitle

\section{\it Introduction}
\label{intro}

Studying objects of a given category we usually wish to have a
notion when these objects are reasonable or appropriate to make
concrete computations. For example, a set can be measurable, a
series can be converging, a vector space can be finite-dimensional,
etc. It is then natural to find such a notion for objects in motivic
categories, for example, in the category of Chow motives $\CHM $ or
in Voevodsky's triangulated category $\DM $ of motives over a field.
Maybe the first step in this direction was made by S.Kimura who has
introduced the notion of finite-dimensionality for Chow motives and
shown its relation to deep problems in algebraic geometry, see
\cite{Kimura}. An algebraic counterpart of Kimura's theory was
independently worked up by P.O'Sullivan, see an overview in
\cite{Andre}.

Y. Andr\'e has shown, \cite{Andre}, \cite{Andre2}, that if a Chow
motive $M$ is finite-dimensional in the sense of Kimura-O'Sullivan,
the corresponding motivic zeta function $\zeta _M(t)$ is rational in
the ring $K_0(\CHM )[[t]]$, where $K_0(\CHM )$ is the Grothen\-dieck
group of the category $\CHM $. In particular, zeta function is
rational for motives of the abelian type in $\CHM $, loc.cit.
Moreover, recently F. Heinloth proved that for those motives zeta
function satisfies a functional equation, \cite{Heinloth}, answering
the question of Andr\'e.

Recall that the notion of motivic zeta function is naturally
connected with many other important topics in different fields of
mathematics. It was introduced by M.Kapranov in \cite{Kapranov} who
proved its rationality and functional equation for smooth projective
curves with respect to any motivic measure $\mu $ with the property
$\mu (\AF ^1)\neq 0$. In particular, if $\mu $ counts the number of
points of a curve defined over a finite field, then the motivic zeta
function is the usual Hasse-Weil zeta function associated with the
curve. The rationality of the Hasse-Weil zeta function for all
varieties was done by Dwork in \cite{Dwork}.

The main goal of the present paper is to show that rationality of
motivic zeta functions can also be considered as a good notion of
measurability for motives in Voevodsky's category $\DM $, possibly
adding functional equations or something. We will prove that
rationality of motivic zeta function has $2$-out-of-$3$ property in
distinguished triangles in the triangulated category $\DMG $ of
geometrical motives over a field. It follows that if $\DMGA $ is a
thick tensor subcategory in $\DMG $ generated by motives of
quasi-projective curves, the motivic zeta function $\zeta _M$ is
rational for any motive $M$ from $\DMGA $. In particular, one can
get an example of a not finite-dimensional motive whose zeta
function is rational. Finally, $2$-out-of-$3$ property of
rationality of motivic zeta function in $\DMG $ leads to a new
conjecture: zeta function $\zeta _M$ is rational for any motive $M$
in the category $\DMG $.

Our approach is homotopical. We use the fact that under mild
assumptions the big category $\DM $ is equivalent to the homotopy
category of an appropriate model monoidal category enriched by
simplicial sets. In that sense it can be considered as an abstract
stable homotopy category whose monoidal product comes from the model
level. Then $2$-out-of-$3$ property for rationality of zeta
functions can be proved using the Postnikoff tower constructed in
\cite{Gu}. Actually, we will build two opposite special
lambda-structures on the Grothendieck group of the thick
triangulated subcategory of compact objects in any homotopy category
of a simplicail model monoidal category localized with $\QQ $.

\section{\it Basics on motivic zeta functions}

Let us start with the folowing important example. Assume that $k=\FF
_q$ is a finite field of $q$ elements, let $X$ be an algebraic
variety over $\FF _q$ and let $\# X(\FF _q)$ be the number of $\FF
_q$-rational points on $X$. Then the usual Hasse-Weil zeta function
$\zeta _X$ associated with $X$ can be defined by the formula
  $$
  \zeta _X(t)=\exp \left( \sum _{n=1}^{\infty }\# X(\FF _{q^n})
  \frac{t^n}{n}\right)
  $$
After the application of the exponential formula coming from
combinatorics,
%\cite[5.1]{Stanley},
it can be rewritten as
  $$
  \zeta _X(t)=\sum _{n=0}^{\infty }\# \Sym ^nX(\FF _q)t^n\; ,
  $$
where $\Sym ^nX$ is the symmetric power of $X$, i.e. the quotient of
$X^{\times n}$ by the action of the symmetric group.

For any field $k$ let $\Var $ be the category of quasi-projective
varieties over $k$. Let, furthermore, $\ZZ [\Var ]$ be the free
abelian group generated by isomorphism classes $[X]$ of
quisi-projective varieties $X$ over $k$. The Grothendieck group
$K_0(\Var )$ of the category $\Var $ is, by definition, the quotient
of $\ZZ [\Var ]$ by the minimal subgroup containing relations of the
type $[X]=[Z]+[X\backslash Z]$ where $Z$ is a closed subvariety in
$X$. Then $K_0(\Var )$ is naturally a commutative ring with unit and
with a product induced by fibered products of varieties over $k$. A
motivic measure $\mu $ is a ring homomorphism $\mu :K_0(\Var )\to A$
to any other commutative ring $A$. Given $\mu $, for any variety $X$
over $k$ we can consider the corresponding zeta-function
  $$
  \zeta _{X,\mu }(t)=\sum _{n=0}^{\infty }\mu [\Sym ^nX]t^n\; ,
  $$
see \cite{Kapranov}.

For example, if $\mu $ counts the number of points over a finite
field, then $\zeta _{X,\mu }(t)$ is the above Hasse-Weil zeta
function of $X$. In that case $\zeta _X(t)$ is rational by Dwork's
result, \cite{Dwork}.

Kapranov proved, see \cite[1.1.9]{Kapranov}, that the motivic zeta
function $\zeta _{X,\mu }(t)$ is rational when $X$ is a smooth
projective curve carrying a divisor of degree one, $A$ is a field
and the motivic measure $\mu (\AF ^1)$ of the affine line is not
zero. On the other hand, Larsen and Lunts have shown in \cite{LL}
that there exists a measure $\mu $ built on the base of Hodge
numbers $h^{i,0}$, such that $\zeta _{X,\mu }$ is not rational for a
surface $X$ with $h^{2,0}\neq 0$.

Now let $\CHM $ be the category of Chow motives over $k$ with
coefficients in $\QQ $. Since $\CHM $ is a symmetric monoidal
additive category, one can construct its Gothendieck ring $K_0(\CHM
)$ in a standard way, i.e. taking direct sums as sums and tensors
products as products in $K_0$. In \cite{GilletSoule} Gillet and
Soul\'e constructed a motivic measure
  $$
  \mu _{GS}:K_0(\Var )\lra K_0(\CHM )
  $$
sending a smooth projective variety $X$ to the class of its Chow
motive $M(X)$. For any given Chow motive $M$ let
  $$
  \zeta _M(t)=\sum _{n=0}^{\infty }[\Sym ^nM]t^n
  $$
be the corresponding zeta function with coefficients $K_0(\CHM )$.
By the result of Del Ba$\tilde {\rm n}$o and Navarro Aznar,
\cite{BNA},
  $$
  \mu _{GS}[\Sym ^nX]=[\Sym ^nM(X)]\; ,
  $$
whence
  $$
  \zeta _{M(X)}=\zeta _{X,\mu _{GS}}
  $$
for any $X$. If $M$ is a Chow motive which is finite-dimensional in
the sense of Kimura-O'Sullivan, then $\zeta _M(t)$ is rational, see
\cite{Andre} and \cite{Andre2}.

A board generalization can be done as follows. Let $\aA $ be any
pseudoabeian symmetric monoidal $\QQ $-linear category with a
monoidal product $\otimes $. Then we have wedge and symmetric powers
of objects $X$ in $\aA $ as images of the corresponding idempotents
in the group algebra of the symmetric group $\sgr _n$ acting on
$X^{\otimes n}$. Let $\ZZ [\aA ]$ be the free abelian group
generated by isomorphisms classes of objects in $\aA $, and let
$K_0(\aA )$ be the Grothendieck group of the category $\aA $, i.e.
the quotient of $\ZZ [\aA ]$ by the minimal subgroup generated by
expressions of type
  $$
  [X\oplus Y]-[X]-[Y]\; .
  $$
Clearly, it has a ring structure induced by the monoidal product in
$\aA $.

Recall that a lambda-structure on a commutative ring $A$ with unit
$1$ is a chain of endomorphisms
  $$
  \lambda ^i:A\to A\; , i=0,1,2,\dots \; ,
  $$
such that $\lambda ^0(a)=1$, $\lambda ^1(b)=r$ and
  $$
  \lambda ^n(a+b)=\sum _{i+j=n}\lambda ^i(a)\lambda ^j(b)
  $$
for all $a$ and $b$ in $A$. It can be also defined as a group
homomorphism
  $$
  \lambda _t:A\lra 1+A[[t]]^+
  $$
from the additive group of $A$ to the multiplicative group
$1+A[[t]]^+$ of formal power series of type $1+a_1t+a_2t^2+\dots $.
The group $1+A[[t]]^+$ has a lambda-structure itself, and the
lambda-structure on $A$ is called to be special if $\lambda _t$
commutes with lambda-operations. Given two lambda structures
$\lambda $ and $\sigma $ on the same ring, they are called to be
opposite if
  $$
  1+\sum _{i=1}^{+\infty }\lambda ^i(a)t^i=
  \left(1+\sum _{i=1}^{+\infty }\sigma ^i(a)(-t)^i\right)^{-1}
  $$
for all $a$ in $A$.

Turning back to the category $\aA $, wedge and symmetric powers in
it give rise to special $\lambda $-structures in the ring $K_0(\aA
)$, opposite each other, \cite[4.1]{Heinloth}. We will denote these
$\lambda $-structures by $\lambda _+$ for wedge and $\lambda _-$ for
symmetric powers respectively. For example, if $X\in \Ob (\aA )$
then
  $$
  \lambda ^n_+[X]=[\wedge ^nX]
  $$
and
  $$
  \lambda ^n_-[X]=[\Sym ^nX]\; .
  $$
Let also
  $$
  \lambda ^{\pm }_t:K_0(\aA )\lra 1+K_0(\aA )[[t]]^+
  $$
be the group homomorphism corresponding to the $\lambda $-structure
$\lambda _{\pm }$. Then for any object $X$ in $\aA $ we can define
its zeta function $\zeta _X(t)$ by the formula:
  $$
  \zeta _X(t):=\lambda ^-_t([X])\; .
  $$
If $\aA $ is the category of Chow motives $\CHM $ then we arrive to
the above motivic zeta function with respect to the measure
constructed by Gillet and Soul\'e.

Below we are mainly interested in motivic measures $\mu $ which can
be factored through $\mu _{GS}$:
  $$
  \mu =\tau \circ \mu _{GS}
  $$
for some $\tau $. In that case, if we know rationality of $\zeta
_{M(X)}$, then we also know rationality of $\zeta _{X,\mu }$.

Any reasonable motivic measure which can be defined in terms of
appropriate cohomology groups can be factored through $K_0(\CHM )$.
For example, given any quasi-projective variety $X$ over $\CC $ we
define its Hodge numbers $h^{p,q}$ as dimensions of the
corresponding bigraded pieces of the mixed Hodge structure on the
cohomology with compact support $H_c^{p+q}(X,\QQ )$ of $X$. Then the
motivic measure sending $X$ to its Hodge polynomial $\sum
h^{p,q}u^pv^q$ can be defined in terms of mixed Hodge realizations.
Therefore it factors through $K_0(\CHM )$. Another interesting
example of a measure factoring through $\mu _{GS}$ can be provided
by conductors of $l$-adic representations over a number field, see
\cite{DL}.

Finally, let us turn to the triangulated setting. Let $\aT $ be a
triangulated category with shift functor
  $$
  X\mapsto \Sigma X\; ,
  $$
and let $\ZZ [\aT ]$ be the free abelian group generated by
isomorphism classes of objects in $\aT $. Let, furthermore, $S(\aT
)$ be the minimal subgroup in $\ZZ [\aT ]$ generated by elements
  $$
  [Y]-[X]-[Z]
  $$
whenever $Y=X\oplus Z$, and let $T(\aT )$ be the minimal subgroup
generated by the same expressions whenever
  $$
  X\lra Y\lra Z\lra \Sigma X
  $$
is a distinguished triangle in $\aT $. The quotient $\ZZ [\aT
]/S(\aT )$ is the above usual Grothendieck group denoted now as
$K_0^{\oplus }(\aT )$, while the quotient $\ZZ [\aT ]/T(\aT )$ is a
more subtle, ``triangulated" Grothendieck group of $\aT $.
Evidently, $K_0(\aT )$ is a quotient of $K_0^{\oplus }(\aT )$ by
$T(\aT )/S(\aT )$.

Note that, if $\aT $ is a derived category of a nice abelian
category $\aA $, then $K_0(\aT )$ is isomorphic to $K_0(\aA )$, so
that the triangulated $K_0$ makes sense.

If we assume that $\aT $ is symmetric monoidal, then $K_0$ of it is
a ring, and the above two subgroups are ideals in $K_0$. Assume,
furthermore, that $\aT $ is $\QQ $-linear and pseudoabelian. Then
wedge and symmetric powers live in the category. Since $\aT $ is
additive, $K_0^{\oplus }(\aT )$ has two canonical lambda-structures
$\lambda _{\pm }$ by Lemma 4.1 in \cite{Heinloth}. A crucial
question is then whether or not these lambda-structures induce
lambda structures on the triangulated group $K_0(\aT )$. The
positive answer to this question would have interesting corollaries
when applying to zeta functions in the triangulated category of
motives over a field.

In the next section we will show an existence of two opposite
special lambda-structures in $K_0(\aT ')$ induced by wedge and
symmetric powers in $\aT $, where $\aT '$ is a thick symmetric
monoidal subcategory generated by compact objects in any abstract
stable homotopy category $\aT $.

\section{\bf \it Main result}
\label{mr}

Let $\aC $ be a pointed simplicial model monoidal category, let co
\cite{Hovey1}, and let
  $$
  \aT =Ho(\aC )
  $$
be the homotopy category of $\aC $. Then $\aT $ is a triangulated
category whose shift functor is induced by the simplicial suspension
  $$
  \Sigma X=X\wedge S^1
  $$
in $\aC $, loc.cit. Such triangulated category has a symmetric
monoidal product $\otimes $ induced by the symmetric monoidal
product $\wedge $ of the category $\aC $, and it can be considered
as an abstract prototype for all reasonable stable homotopy
categories, see \cite{Hovey1}.

Since we are interested in the study of wedge and symmetric powers
of objects in $\aT $, we will assume that $\aT $ is $\QQ
$-localized. If we look on objects in $\aC $ as on spectra, then we
consider the homotopy category of $\QQ $-local spectra.

There are several examples of such triangulated categories arising
in algebraic topology and motivic algebraic geometry. The homotopy
category of $\QQ $-local topological symmetric spectra over a point
is just the category of graded $\QQ $-vector spaces. This is not
interesting. However, the rational stable homotopy theory of
$S^1$-equivariant spectra is still interesting, see
\cite{Greenlees}.

But the main example for our purposes is the Morel-Voevodsky
homotopy category $\SH (X)$ of $\QQ $-local motivic symmetric
spectra over a Noetherian base scheme $X$, see \cite{Jardine} and
\cite{VoevodskyBerlin}. If $X=\Spec (k)$ we will write $\SH $
instead of $\SH (\Spec (X))$.

As it was announced by Morel, \cite{Morel}, if $-1$ is a sum of
squares in the ground field $k$, the category $\SH $ is equivalent
to the big category $\DM $ of triangulated motives over $k$.

Let $H_{\ZZ }$ be the motivic Eilenberg-MacLane spectrum inducing
motivic cohomology of schemes over $k$. By the result in \cite{RO},
the category $\DM $ is also equivalent to the homotopy category
$Ho(H_{\ZZ }-\modules )$ of modules over the ring spectrum $H_{\ZZ
}$ when $k$ is of characteristic zero.

Using these equivalences we can apply any result obtained in $\aT
=Ho(\aC )$ to the category $\DM $.

Since $\aT $ is a homotopy category, it has all coproducts. Then
$K_0(\aT )=0$ by Eilenberg swindle. Indeed, let $X$ be an object in
$\aT $ and let $[X]$ be its class in $K_0$. Let
  $$
  Y=X\coprod X\coprod \dots
  $$
be the coproduct of a countable number of copies of $X$ in $\aT $.
Then
  $$
  X\oplus Y=Y\; ,
  $$
whence $[X]=0$ in $K_0$. Therefore, dealing with $\aT =Ho(\aC )$ it
is reasonable to work with the thick triangulated subcategory
  $$
  \aT ':=\aT ^{\aleph _0}
  $$
of compact objects in $\aT $, see \cite{Neeman}. For example, if
$\aT =\DM $ over a field then $\aT '=\DM '$ is nothing but the
triangulated category of geometrical motives $DM_{gm}$ over $k$, see
\cite{VoevodskyTriang}.

A result is then as follows:

\begin{theorem}
\label{main} Wedge and symmetric powers induce two special
lambda-structures in the ring $K_0(\aT ')$, which are opposite each
other.
\end{theorem}

Essentially, the below proof of this theorem is based on the
existence of a Postnikoff tower connecting wedge (symmetric) powers
of the vertex $Y$ in a given distinguished triangle
  $$
  X\to Y\to Y\to \Sigma X
  $$
with wedge (symmetric) powers of two another vertices $X$ and $Z$.
Without loss of generality, applying cofibrant replacement, we can
assume that both $X$ and $Y$ are cofibrant and the above
distinguished triangle is a cofibration triangle, so that $Z=Y/X$
may be considered as a ``contraction of $X$ inside $Y$ into a
point".

Let $n$ be a natural number and let $K$ be an $n$-cube commutative
diagram in $\aC $ whose vertices are indexed by $n$-tuples $v=\{
v_1,\dots ,v_n\} $, $v_i=0,1$, and built up as smash-products
$A_1\wedge \dots \wedge A_n$ in $\aC $ with $A_j=X$ if $v_j=0$ and
$A_j=Y$ if $v_j=1$. The morphisms of $K$ are smash products of
finite collections of copies of the morphism $X\to Y$ and the
identity $X\to X$, according to vertices in $K$. For example, when
$n=2$ the diagram $K$ looks like
  $$
  \diagram
  X\wedge X \ar[dd]^-{f\wedge 1}
  \ar[rr]^-{1\wedge f} & &
  X\wedge Y \ar[dd]^-{f\wedge 1} \\ \\
  Y\wedge X \ar[rr]^-{1\wedge f} & & Y\wedge Y
  \enddiagram
  $$
where the objects $X\wedge Y$ and $Y\wedge X$ correspond to the
vertices $(0,1)$ and $(1,0)$ respectively.

For any $0\leq i\leq n$ let $K^i$ be the subdiagram in $K$ generated
by vertices indexed by $n$-tuples $v$ containing $\leq i$ of units.
Then we have the filtration
  $$
  K^0\subset K^1\subset \dots \subset K^n=K\; .
  $$
Let
  $$
  D^i=\colim K^i
  $$
be the colimit of the diagram $K^i$ in the category $\aC $. Clearly,
$D^0=X^{\otimes n}$ and $D^n=Y^{\otimes n}$. For each $i$ the
inclusion $K^i\subset K^{i+1}$ induces a morphism
  $$
  w_i:D^i\lra D^{i+1}
  $$
on the colimits in $\aC $. Let
  $$
  E^{i+1}=\cone (w_i)\; .
  $$
Any permutation $\sigma \in \sgr _n$ gives rise to an endomorphism
$\Gamma _{\sigma }$ of the distinguished triangle
  $$
  D^i\stackrel{w_i}{\lra }D^{i+1}\lra E^{i+1}\lra \Sigma D^i\; .
  $$
These $\Gamma _{\sigma }$ are induced by the permutation $\sigma $
on colimits. Considering the idempotents of the group algebra $\QQ
\sgr _n$ inducing wedge (symmetric) powers in $\aT $, we get
distinguished triangles
  $$
  \diagram
  I_+^i \ar[r]^-{} &
  I_+^{i+1} \ar[r]^-{} &
  J_+^{i+1} \ar[r]^-{} & \sgr I_+^i
  \enddiagram
  $$
for the alternated case, and a similar distinguished triangles
  $$
  \diagram
  I_-^i \ar[r]^-{} &
  I_-^{i+1} \ar[r]^-{} &
  J_-^{i+1} \ar[r]^-{} & \sgr I_-^i
  \enddiagram
  $$
in the symmetric case. It is easy to see that
  $$
  I_+^0=\wedge ^nX\; ,\; \; \; I_+^n=\wedge ^nY\; ,
  $$
and, similarly,
  $$
  I_-^0=\Sym ^nX\; ,\; \; \; I_-^n=\Sym ^nY\; .
  $$
The key point is that the above cones $C_{\pm }^i$ can be computed
by the rule:
  $$
  J_+^i=\wedge ^iZ\otimes \wedge ^{n-i}X
  $$
and
  $$
  J_-^i=\Sym ^iZ\otimes \Sym ^{n-i}X\; .
  $$
The detailed proof of this fact is given in \cite{Gu}.
%Let
%  $$
%  P_{\pm }=\{ P^i_{\pm }\} _{i=0,1,\dots ,n}
%  $$
%be the above Postnikoff towers.

Now we need a few elementary but useful algebraic observations. Let
$M$ be an abelian monoid and let $M^+$ be the group completion of
$M$. The canonical morphism $M\to M^+$, $m\mapsto [m]$ is universal
with respect to all morphisms from $M$ to abelian groups inverting
non-zero elements from $M$.

Assume we are given with an equivalence relation $R$ on $M$. We will
say that $R$ is additive if the following condition holds: for any
two $a$, $a'$, $b$ and $b'$ in $M$, if $(a,a')\in R$ and $(b,b')\in
R$, then it follows that $(a+a',b+b')$ is also in $R$. If $R$ is
additive then we can construct a quotient additive monoid $M/R$.

Let $R$ be an additive equivalence relation on $M$ and let $R^+$ be
a subgroup in $M$ generated by elements $[a]-[a']$, such that
$(a,a')\in R$. Then the abelian group $M^+/R^+$ is canonically
isomorphic to the abelian group $(M/R)^+$ because both compositions
$M\to M^+\to M^+/R^+$ and $M\to M/R\to (M/R)^+$ are universal with
respect to morphisms $f:M\to A$ of the monoid $M$ to abelian groups
$A$, such that $f(a)=f(a')$ as soon as $(a,a')\in R$.

Let $S$ be a set and let $\NN [S]$ be a free abelian monoid
generated by $S$. Let $\rho $ be a reflexive and symmetric relation
on $\NN [S]$. If $(a,b)\in \rho $ the we will say that $a$ is
elementary $\rho $-equivalent to $b$. Build an additive equivalence
relation $\langle \rho \rangle $ on $\NN [S]$ generated by $\rho $
as follows: two linear combinations $a$ and $a'$ from $\NN [S]$ are
called to be equivalent if there exist a chain of linear
combinations $a_0,a_1,\dots a_n$, such that $a_0=a$, $a_n=a'$ and
for each $i$ the element $a_{i+1}$ can be obtained from $a_i$ by a
replacement of a summand in $a_i$ by a $\rho $-equivalent summand.
For short, let $\NN [S]/\rho $ be the quotient of $\NN [S]$ by
$\langle \rho \rangle $ and $\rho ^+=\langle \rho \rangle ^+$.
Certainly, $(\NN [S])^+=\ZZ [S]$. From the previous observation we
have:
   $$
   (\NN [S]/\rho )^+=\ZZ [S]/\rho ^+\; .
   $$

For example, let $\aA $ be an additive category, $\NN [\aA ]:=\NN
[\Ob (\aA )]$ and let $\sigma $ be the set of pairs $([X\oplus
Y],[X]+[Y])$ and their transposes. We see that $\ZZ [\aA ]/\sigma
^+$ is the additive Grothendieck group $K_0^{\oplus }(\aA )$ of $\aA
$. The above isomorphism shows then that $K_0^{\oplus }(\aA )$ can
be also described as a completion of the monoid $\NN [\aA ]/\sigma
$.

In the triangulated situation we have the following. Let $\aT $ be a
triangulated category, $\NN [\aT ]:=\NN [\Ob (\aT )]$ and let
$\Delta $ be the set of pairs $([Y],[X]+[Z])$, where $XYZ$ is a
distinguished triangle, and their transposes. Again we see that $\ZZ
[\aT ]/\Delta ^+$ is the triangulated Grothendieck group $K_0(\aT )$
and the above isomorphism shows then that it can be also described
as a completion of the monoid $\NN [\aT ]/\Delta $.

In order to finish the proof of Theorem \ref{main} we will consider
the case of wedge powers only. The symmetric case is analogous.

For any object $X$ in $\aT $ let $[X]$ be its class in the
Grothendieck group, $K_0^{\oplus }$ or $K_0$. We have a countable
set of maps
  $$
  \lambda ^n:K_0^{\oplus }(\aT )\lra K_0^{\oplus }(\aT )
  \; ,\; \; n=0,1,2,\dots \; ,
  $$
such that
  $$
  \lambda ^n[X]=[\wedge ^nX]
  $$
for any $X$ in $\aT $, and the value of $\lambda ^n$ on a sum
$[X]+[Y]$ is determined by the K\"unneth rule. These are the
$\lambda $-operations considered in \cite{Heinloth}. In order to
define $\lambda $-operations by the same rule in $K_0(\aT )$ we need
only to show that the above maps $\lambda ^i$ respect also the
subgroup $T(\aT )$. Since $K_0(\aT )=(\NN [\aT ]/\Delta )^+$, it is
enough to show that if two linear combinations $a$ and $b$ are
elementary $\Delta $-equivalent, then $\lambda ^na$ is $\Delta
$-equivalent to $\lambda ^nb$.

Without loss of generality we can assume that $a=[Y]$, $b=[X]+[Z]$
and we have a distinguished triangle
  $$
  X\lra Y\lra Z\lra \Sigma Z
  $$
in $\aT $. For each $n$ let
  $$
  \diagram
  I_+^i \ar[r]^-{} &
  I_+^{i+1} \ar[r]^-{} &
  J_+^{i+1} \ar[r]^-{} & \sgr I_+^i\; ,
  \enddiagram
  $$
  $$
  i=0,1,\dots ,n-1\; ,
  $$
be the Postnikoff system as above. Let
  $$
  a_i=[I_+^i]\; ,\; \; \; c_i=[J_+^i]\; ,\; \; \; x_i=[\wedge ^iX]\; ,
  \; \; y_i=[\wedge ^iY]\; \; \; \hbox{and}\; \; \; z_i=[\wedge ^iZ]
  $$
for all non-negative integers $i$. The above Postnikoff tower can be
rewritten in terms of the ring $K_0(\aT )$ as follows:
  $$
  a_0=x_n\; ,\; \; \; a_n=y_n\; ,
  $$
  $$
  c_i=z_ix_{n-i}\; ,\; \; i=0,1,\dots ,n\; ,
  $$
and
  $$
  a_{i+1}=a_i+c_{i+1}\; ,\; \; i=0,1,\dots ,n-1\; ,
  $$
whence we get:
  $$
  y_n=\sum _{i=0}^nz_ix_{n-i}
  $$
in the ring $K_0(\aT )$, where $y_n=\lambda ^n[Y]$, $x_0=1$ and
$z_0=1$. It means that the element $\lambda ^n([X]+[Z])$ is $\Delta
$-equivalent to the element $\lambda ^n[Y]$ in the ring $K_0^{\oplus
}$, as desired.

%Since the lambda-operations on $K_0$ are induced by the
%lambda-operations on $K_0^{\oplus }$, we have the commutative square
%  $$
%  \diagram
%  K_0^{\oplus }(\aT ') \ar[dd]^-{}
%  \ar[rr]^-{\lambda ^i} & &
%  K_0^{\oplus }(\aT ') \ar[dd]^-{} \\ \\
%  K_0(\aT ')\ar[rr]^-{\lambda ^i} & &
%  K_0(\aT ')
%  \enddiagram
%  $$
%for each index $i$.

Let $\Lambda ^i$ be the $i$-th operation of the canonical
lambda-structure on $1+K_0(\aT ')[[t]]]^+$. The precompositions of
the two paths $\Lambda ^i\lambda _t$ and $\lambda _t\lambda ^i$ in
the diagram
  $$
  \diagram
  K_0(\aT ') \ar[dd]^-{\lambda _t}
  \ar[rr]^-{\lambda ^i} & &
  K_0(\aT ') \ar[dd]^-{\lambda _t} \\ \\
  1+K_0(\aT ')[[t]]^+ \ar[rr]^-{\Lambda ^i} & &
  1+K_0(\aT ')[[t]]^+
  \enddiagram
  $$
with the surjective homomorphism
  $$
  K_0^{\oplus }(\aT ')\lra K_0(\aT ')
  $$
coincide because the $\lambda $-structure on $K_0^{\oplus }$ is
special. Therefore the diagram is commutative, i.e. the induced
$\lambda $-structure on the $K_0$-level is special as well.

By the same argument, the lambda-structure on $K_0$ defined by wedge
powers is opposite to the $\lambda $-structure defined by symmetric
powers.

\section{\it Some applications}

A formal power series $\xi (t)$ in variable $t$ with coefficients in
a commutative ring $A$ is called rational, if there exists two
polynomials $a(t)$ and $b(t)$ in $A[t]$, such that $a\xi =b$ in
$A[[t]]$.

For any element $a\in K_0(\aT ')$ let $\zeta _a(t)=\lambda ^-_t(a)$
be the corresponding zeta function of $a$. If  $a=[X]$ for some
object $X$ in $\aT '$, then we write $\zeta _X$ instead of $\zeta
_a$. The suspension $\Sigma :\aT \to \aT $ induces an involution
$a\mapsto -a$ on $K_0(\aT ')$. Theorem \ref{main} gives $\zeta
_{-a}=\zeta _a^{-1}$ for any element $a$ in $K_0$. In particular,
$\zeta _{\Sigma X}=(\zeta _X)^{-1}$ for any object $X$ in $\aT '$.
It follows that the suspension does not change rationality of the
corresponding zeta function. Also we have:

\begin{corollary}
\label{23} Let $X\to Y\to Z\to \Sigma X$ be a distinguished triangle
in $\aT $. If two out of three zetas $\zeta _X$, $\zeta _Y$ and
$\zeta _Z$ are rational functions, then the third zeta function is
also rational.
\end{corollary}

\begin{pf}
Since $XYZ$ is a distinguished triangle, $[Y]=[X]+[Z]$ in $K_0(\aT
)$. Since we work with zeta-functions induced by $\lambda
$-structures in $K_0$,
  $$
  \zeta _Y=\zeta _X\cdot \zeta _Z\; ,
  $$
whence the proof.
\end{pf}

Assume now that either $\cha (k)=0$ or $-1$ is a sum of squares in
$k$. In that case, as we have seen in Section \ref{mr}, the category
$\DM $ can be interpreted as a category $\aT =Ho(\aC )$ via the
equivalence with either the motivic stable category $\SH $ or with
the homotopy category of modules over the motivic Eilenberg-MacLane
spectrum. Therefore, Theorem \ref{main} can be applied to the
category $\DM '$ as well. For any motive $M\in \DM '$ let $ \zeta
_M(t)$ be the corresponding motivic zeta function of $M$.

The advantage of Theorem \ref{main} is that it extends the range of
varieties whose motivic zeta function is rational. Indeed, let $\DM
_{\leq 1}^{\otimes }$ be a thick symmetric monoidal subcategory in
$\DM $ generated by motives of quasi-projective curves over $k$.

\begin{corollary}
\label{one} For any motive $M$ in $\DMGA $ its zeta function $\zeta
_M(t)$ is rational.
\end{corollary}

\begin{pf}
The motive $M(X)$ of a quasi-projective curve $X/k$ is
finite-dimensional in the sense of Kimura-O'Sullivan, see \cite{Gu}
or \cite{Mazza}. Therefore, $\zeta _{M(X)}$ is strictly rational by
Andr\'e's result, see \cite[13.3.3.1]{Andre2}. If $M$ is in $\DM
_{\leq 1}^{\otimes }$, then its class $[M]$ belongs to a subring in
$K_0(\DM ')$ generated by classes of motives of quasi-projective
curves.

Given a $\lambda $-structure on a ring $A$, if $\zeta _a$ and $\zeta
_b$ are rational for two elements $a$ and $b$ in $A$, then $\zeta
_{a+b}$ is also rational. If we assume, furthermore, that this
lambda-structure is special, then $\zeta _{ab}=\zeta _a\star \zeta
_b$ where $\star $ is the product in the lambda-ring $1+A[[t]]^+$.
At the same time, for any two elements $f(t)$ and $g(t)$ in
$1+A[[t]]^+$ if $f$ and $g$ are ratio of polynomials then so is the
product $f\star g$.

Since the lambda-structure on $K_0(\DM ')$ is special, rationality
of zeta functions is closed under sums and products of elements in
$K_0(\DM ')$, whence the result.
\end{pf}

Corollary \ref{23} says that motivic zeta function has
$2$-out-of-$3$ property in distinguished triangles in $\DM '$. This
correlates with Lemma 3.1 in \cite{LL} via Gyzin distinguished
triangles in the category $\DM '$. In general, there are several
canonical distinguished triangles in $\DM $ each of which gives new
examples of varieties whose motivic zeta function is rational.

\begin{example}
{\sf Let $X$ be a K3 surface over $\CC $, such that $M(X)$ is in
$\DMGA $. This is so if $X$ is, for example, of Kummer type or a
Weil quartic in $\PR ^3$. Then $\zeta _X$ is rational by Corollary
\ref{one}. As it was shown by O'Sullivan, \cite{Mazza}, there exists
a Zariski open $U$ in $X$, such that the motive $M(U)$ is not
finite-dimensional. However, the function $\zeta _U$ is rational.
Indeed, the complement $Z=X-U$ is a union of curves. The motive of
any quasi-projective curve is finite-dimensional, whence rationality
of $\zeta _Z$. Since $\zeta _X$ is rational, applying Corollary
\ref{23} to the Gyzin distinguished triangle in $\DM $ associated
with the pair $X,U$ we see that $\zeta _U$ is rational.}
\end{example}

\begin{example}
{\sf Let $X=X_1\cup \dots \cup X_n$ be a union of quasi-projective
surfaces whose zeta functions $\zeta _{X_i}$ are rational. Then
$\zeta _X$ is rational. Indeed, if $n=2$, we apply Corollary
\ref{23} and rationality of zeta functions of quasiprojective curves
to the corresponding Mayer-Vietoris distinguished triangle in $\DM
$. If $n>2$, we just write down an appropriate Postnikoff tower and
apply Corollary \ref{23} several additional times.}
\end{example}

\begin{example}
{\sf Let $X$ be a quasi-projective variety over $k$, let $Z$ be a
closed subvariety in $X$ and let $f:X'\to X$ be a blow up of $X$ in
$Z$. Assume that $\zeta _Z$ is rational. Then we use the
corresponding blow up distinguished triangle in $\DM $ in order to
show that $\zeta _X$ is rational if and only if $\zeta _{X'}$ is
rational. Of course, if all appearing varieties are smooth
projective, then the same claim follows from the result in
\cite{Heinloth} and Manin's motivic formula for blow ups in the
category of Chow motives over $k$.}
\end{example}

\begin{example}
{\sf Let $C$ be a smooth projective curve, let $Y$ be a quotient of
$C\times C\times C$ by the cyclic action of $\ZZ /3$, and let $X$ be
resolution of singularities of $Y$. Then $X$ is a smooth projective
threefold whose motive is in $\DMGA $ and the motivic $\zeta _X$ is
rational.}
\end{example}

\begin{remark}
{\sf By recent result of Bondarko, \cite{Bondarko}, $K_0(\CHM )$ is
isomorphic to the triangulated $K_0(\DM ')$ if we consider both
categories over a field of characteristic zero. Most probably this
isomorphism also gives the two above canonical lambda-structures in
Voevodsky's category $\DM '$. However, Theorem \ref{main} allows to
work with lambda-structures in the 2-functor $\SH (-)$ on schemes in
general. So, on the one hand, we can use $\SH (-)$ in order to build
new varieties whose motivic zeta functions are rational, on the
other hand, we can use Bondarko's isomorphism in order to factor
interesting motivic measures through $K_0(\DM ')$.}
\end{remark}

\begin{example}
{\sf Assume we are given with a nice motivic measure $\mu $ which
factors through $K_0(\CHM )$, say $\mu $ is made by Hodge
polynomials of varieties over $\CC $ or by conductors of $l$-adic
representations for \'etale cohomology of varieties over $\QQ $.
Since $K_0(\CHM )$ is naturally isomorphic to $K_0(\DM ')$, we can
also factor $\mu $ through $K_0(\DM ')$ and apply Corollary
\ref{one}. Then we have that $\zeta _{X,\mu }$ is rational for any
variety $X$, such that $M(X)$ is in $\DMGA $.}
\end{example}

\begin{remark}
\label{Ayoub} {\sf Since $\DM '$ is generated by $\CHM $, Theorem
\ref{main} shows that if $\zeta _M$ is rational for all Chow motives
over $k$, then $\zeta _M$ is rational for all motives in $\DM '$.
Let $\eta $ be the generic point of the affine line $\AF ^1$, and
let $0$ be a closed point on it. Let
  $$
  \Psi :\DM '(\eta )\lra \DM '(0)
  $$
be the nearby cycle motivic functor constructed by J.Ayoub in
\cite{Ayoub}. The category $\DM '$ is generated by motives of type
$\Psi M(X)$ where $X$ is a smooth hypersurface in a projective space
over $\eta $, loc.cit. Therefore, in order to show rationality of
motivic zetas for all motives in $\DM '$ it is enough to show that
$\zeta _X$ is rational for a smooth projective hypersurface $X$ in
$\PR ^n$ over $\eta $. The standard approach is as follows. Let
$\mathscr X\to \PR ^N$ be a universal family of smooth hypersurfaces
of degree $d$ over a field. Draw a line $\AF ^1$ in the parameter
space $\PR ^N$ through a smooth hypersurface $X$ and a Weil type
hypersurface $W$. Restricting the universal family on $\AF ^1$ we
get a smaller family $\mathscr Y\to \AF ^1$ which deforms $X$ into
$W$. The motivic $\zeta _W$ is rational, of course. Then one can try
to deform rationality of $\zeta _W$ to rationality of $\zeta _X$
using the homomorphism
  $$
  K_0\Psi :K_0(\DM '(\eta ))\lra K_0(\DM '(0))\; .
  $$
This is, actually, nothing else but a $K_0$-simplification of
Ayoub's approach to Schur-finiteness of $\DM '$. Note that if $M$ is
a Chow motive over $\eta $, then we can use the nice properties of
the functor $\Psi $ established in \cite{Ayoub} in order to compute
$K_0\Psi [M]$ in terms of the usual Chow specialization of the
motive $M$ to $0$.}
\end{remark}

\medskip

{\large Acknowledgements.} I thank Joseph Ayoub for several useful
discussions along the theme of this paper. The work was supported by
NSF Grant DMS-0111298. I am grateful to the Institute for Advanced
Study at Princeton for the support and hospitality in 2005-2006
academical year.

\medskip

\begin{small}

\end{small}

\vspace{5mm}
%\bigskip

{\footnotesize

% \noindent {\sc School of Mathematics, The Institute for Advanced
% Study, Einstein Drive, Princeton, NJ 08540}

\noindent {\it E-mail}: {guletski@ias.edu}

}

\end{document}